\documentclass[12pt]{amsart}
\usepackage{amssymb}
\usepackage{amsbsy}
\usepackage{amscd}
\advance\textheight by \topskip
\makeatletter
\def\cal{\mathcal}
\def\Bbb{\mathbb}

\newenvironment{pf*}[1]{\proof[#1]}{\endproof}
\newcommand{\rom}{\textup}
\renewcommand{\thesubsection}{\thesection(\@roman\c@subsection)}
\makeatother
\newtheorem{Theorem}[equation]{Theorem}
\newtheorem{Corollary}[equation]{Corollary}
\newtheorem{Lemma}[equation]{Lemma}
\newtheorem{Proposition}[equation]{Proposition}

\theoremstyle{definition}
\newtheorem{Definition}[equation]{Definition}

\makeatletter
\renewcommand\section{\@startsection{section}{1}%
  {\z@}{.7\linespacing\@plus\linespacing}{.5\linespacing}%
  {\reset@font\normalfont\bfseries\centering}}
\makeatother

\theoremstyle{remark}
\newtheorem{Remark}[equation]{Remark}

\newtheorem*{Acknowledgments}{Acknowledgments}
\numberwithin{equation}{section}
\numberwithin{figure}{section}

\newcommand{\thmref}[1]{Theorem~\ref{#1}}
\newcommand{\secref}[1]{\S\ref{#1}}
\newcommand{\lemref}[1]{Lemma~\ref{#1}}
\newcommand{\propref}[1]{Proposition~\ref{#1}}
\newcommand{\remref}[1]{Remark~\ref{#1}}

\newcommand{\subsecref}[1]{\S\ref{#1}}

\newcommand{\Romnum}[1]{\expandafter\uppercase\expandafter{\romannumeral #1}} 
\newcommand{\C}{{\Bbb C}}
\newcommand{\Z}{{\Bbb Z}}

\newcommand{\R}{{\Bbb R}}

\newcommand{\CP}{\operatorname{\C P}}

\newcommand{\SO}{\operatorname{\rm SO}}
\newcommand{\U}{\operatorname{\rm U}}
\newcommand{\Pin}{\operatorname{\rm Pin}}
\newcommand{\Spin}{\operatorname{\rm Spin}}
\newcommand{\Spinc}{\Spin^{c}}
\newcommand{\STop}{\operatorname{\rm STop}}
\newcommand{\SpinTop}{\operatorname{\rm SpinTop}}

\newcommand{\Map}{\operatorname{Map}}
\newcommand{\map}{\operatorname{map}}
\newcommand{\Hom}{\operatorname{Hom}}

\newcommand{\Iso}{\operatorname{Iso}}

\newcommand{\im}{\mathop{\text{\rm im}}\nolimits}
\newcommand{\rank}{\operatorname{rank}}
\newcommand{\Sign}{\operatorname{Sign}}
\newcommand{\rel}{\operatorname{rel }}
\newcommand{\celldim}{\operatorname{cell-dim}}

\newcommand{\ind}{\mathop{\text{\rm ind}}\nolimits}

\newcommand{\colim}{\operatorname{colim}}
\newcommand{\G}{{\cal G}} 
\newcommand{\CC}{{\cal C}} 
\newcommand{\UU}{{\cal U}} 
\newcommand{\OO}{{\cal O}} %
\newcommand{\D}{{\cal D}} %
\newcommand{\QQ}{{\cal Q}} %

\newcommand{\id}{\operatorname{id}}

\newcommand{\fb}{\operatorname{BF}}

\begin{document}
\title[Bauer-Furuta invariants under $\Z_2$-actions]{Bauer-Furuta invariants under $\Z_2$-actions}
\author{Nobuhiro Nakamura}
\address{Graduate School of Mathematical Sciences, University of Tokyo, 3-8-1, Komaba, Meguro-ku, Tokyo, 153-8914, Japan}
\email{nobuhiro@ms.u-tokyo.ac.jp}
\begin{abstract}
S.~Bauer and M.~Furuta defined a stable cohomotopy refinement of the Seiberg-Witten invariants.
In this paper, we prove a vanishing theorem of Bauer-Furuta invariants for $4$-manifolds with smooth $\Z_2$-actions. 
As an application, we give a constraint on smooth $\Z_2$-actions on homotopy $K3\#K3$, and construct a nonsmoothable locally linear $\Z_2$-action on $K3\#K3$.
We also construct a nonsmoothable locally linear $\Z_2$-action on $K3$.
\end{abstract}
%
%
%
%
\maketitle
%
\section{Introduction}\label{sec:intro}
%
%
F.~Fang proved a mod $p$ vanishing theorem of Seiberg-Witten invariants under a cyclic group action of prime order $p$~\cite{Fang}.
This theorem, later extended by the author~\cite{Np}, is used when X.~Liu and the author constructed  nonsmoothable locally linear group actions on elliptic surfaces in \cite{LN,LN2}.

On the other hand, S.~Bauer and M.~Furuta defined a stable cohomotopy refinement of the Seiberg-Witten invariants~\cite{BF, B1}.
In general, Bauer-Furuta invariants may be non-trivial even if ordinary Seiberg-Witten invariants are trivial~\cite{B0, B1, FKM}.

In this paper, we investigate Bauer-Furuta invariants under involutions.
In fact, we prove a vanishing theorem of Bauer-Furuta invariants for $4$-manifolds with smooth $\Z_2$-actions.
As an application, we give a constraint on smooth $\Z_2$-actions on homotopy $K3\#K3$, and construct an example of locally linear action on $K3\#K3$ which can not be smoothed. 

To state out results, we need some preliminaries.
Let $X$ be an oriented  smooth $4$-manifold with an orientation-preserving smooth $\Z_2$-action. 
Fixing a $\Z_2$-invariant metric on $X$, we have a $\Z_2$-action on the frame bundle.
Suppose that the $\Z_2$-action on $X$ lifts to a $\Spinc$-structure $c$ over $X$.
Fix a $\Z_2$-invariant connection $A_0$ on the determinant line bundle $L$ of $c$. 
Then the Dirac operator $D_{A_0}$ associated to $A_0$ is $\Z_2$-equivariant, and the $\Z_2$-index of $D_{A_0}$ can be written as $\ind_{\Z_2} D_{A_0} = k_+\C_+ + k_-\C_- \in R(\Z_2)\cong\Z [t]/(t^2-1)$, where $\C_+$ (resp. $\C_-$) is the complex $1$-dimensional representation on which the generator of $\Z_2$ acts by multiplication of $+1$ (resp. $-1$), and $R(\Z_2)$ is the representation ring of $\Z_2$.

Let $b_i$ be the $i$-th Betti number of $X$ and $b_+$  the rank of a maximal positive definite subspace $H^+(X;\R)$ of $H^2 (X;\R)$. 
For any group $G$ and any $G$-space $V$, let $V^G$ be the fixed point set of the $G$-action. Let $b_+^G = \dim H^+(X;\R)^G$. 

Our main theorem is:
\begin{Theorem}\label{thm:main}
Suppose the following conditions are satisfied{\rom :}
\begin{enumerate}
\item $b_1=0$, $b_+\geq 2$ and $b^{\Z_2}_+\geq 1$,
\item $b_+-b_+^{\Z_2}$ is odd,
\item $d(c)=2(k_+ + k_-) - (1 + b_+)=1$,
\item $2k_{\pm}< 1+ b_+^{\Z_2}$.
\end{enumerate}
Then the Bauer-Furuta invariant of $c$ is zero.
\end{Theorem}
\begin{Remark}
The number $d(c)$ is the virtual dimension of the Seiberg-Witten moduli space for the $\Spinc$-structure $c$. 
The Bauer-Furuta invariant is  $\Z/2$-valued if $d(c)=1$ and $k_++k_-$ is even, and is always $0$ if $d(c)=1$ and $k_++k_-$ is odd~\cite{BF}.   
Note also that, when $d(c)=1$ and $b_1=0$, the ordinary Seiberg-Witten invariant is trivial by definition.
\end{Remark}

The strategy of the proof of \thmref{thm:main} is inspired by Bauer's arguments in \cite{B2}. 
The monopole map of a spin structure is $\Pin(2)$-equivariant, and to analyse the monopole map,  Bauer made use of the equivariant obstruction theory on Bredon cohomology. 
In fact, he proved that images of the forgetting map from a $\Pin(2)$-equivariant stable cohomotopy group to  $S^1$-equivariant one are multiples of $2$, which implies that the value of the Seiberg-Witten invariant is even in some situations.

When a $\Z_2$-action on the $\Spinc$-structure is given, we have a $\Z_2\times S^1$-action on the monopole map, instead of the above $\Pin(2)$-action.
By calculating Bredon cohomology groups explicitly, we can show the required vanishing of the Bauer-Furuta invariant.

As an application of \thmref{thm:main}, we construct a nonsmoothable locally linear $\Z_2$-action on $K3\#K3$. ({\it Cf.} \cite{LN,LN2}.)
\begin{Theorem}\label{thm:nonsmoothable}
There exists a locally linear $\Z_2$-action on $X=K3\#K3$ which can not be smooth with respect to any smooth structure on $X$.
\end{Theorem}
To prove \thmref{thm:nonsmoothable}, we basically follow the strategy of \cite{LN,LN2};
In the first step, we give a constraint on smooth actions.
In the second step, we construct a locally linear action which would violate the constraint. 

To obtain a constraint on smooth actions, we use \thmref{thm:main} together with the non-vanishing result of the Bauer-Furuta invariant of $K3\#K3$ by Furuta, Kametani and Minami~\cite{FKM}.
On the other hand, to construct a locally linear action, we invoke a realization theorem due to A.~L.~Edmonds and J.~H.~Ewing~\cite{EE}.

As a byproduct of our argument, we also have:
\begin{Theorem}\label{thm:K3}
There exists a locally linear $\Z_2$-action on $K3$ with isolated fixed points satisfying $b_+^{\Z_2}=3$ which can not be smooth with respect to any smooth structure on $K3$.
\end{Theorem}
\begin{Remark}
J.~Bryan proved every smooth  $\Z_2$-action on $K3$ with isolated fixed points satisfies $b_+^{\Z_2}=3$~\cite{Bryan}. 
Moreover, the nonsmoothable involution in \thmref{thm:K3} has the same fixed point data and action on the $K3$ lattice as a symplectic automorphism of order $2$ on $K3$ called Nikulin involution(\cite{Morrison}, Section 5).
\end{Remark}
The proof of \thmref{thm:K3} is parallel to that of \thmref{thm:nonsmoothable}. 
The difference is that the $G$-spin theorem gives a sufficient constraint on smooth involutions in this case.

The paper is organized as follows.
In Section $2$, we give a brief review on Bauer-Furuta invariants, in particular, a description of the invariants as obstruction classes.
In Section $3$, we review on the equivariant obstruction theory.
In Section $4$, we prove our main theorem (\thmref{thm:main}). 
In Section $5$, we explain the applications (\thmref{thm:nonsmoothable} and \thmref{thm:K3}).
\begin{Acknowledgments}
The author would like to express his gratitude to M.~Furuta, Y.~Kametani and H.~Sasahira for helpful discussions during the completion of this work.
It is also pleasure to thank referees for many valuable suggestions, in particular, for pointing out that \thmref{thm:K3} can be proved without gauge theory.
\end{Acknowledgments}
%
%
%
\section{Bauer-Furuta invariants}\label{sec:bf}
%
%
The purpose of this section is to give a brief review on Bauer-Furuta invariants. (See \cite{BF,B1,Sz} for details.)
%
%
\subsection{Equivariant Bauer-Furuta invariants}\label{subsec:eqbf}
%
%
Suppose that a $4$-manifold $X$ with a $\Z_2$-action satisfies conditions in \thmref{thm:main}.
Let $S^+$ and $S^-$ be the positive and negative spinor bundle of the $\Spinc$-structure $c$, and $L$  the determinant line bundle.

The Seiberg-Witten equations are a system of equations for $\U(1)$-connections $A$ on $L$ and positive spinors $\phi\in\Gamma (S^+)$,
\begin{equation}\label{eq:SW}
\left\{
\begin{gathered}
 D_A \phi = 0,\\
 F_A^+  = q(\phi),
\end{gathered}
\right.
\end{equation}
where $D_A$ denotes the Dirac operator, $F_A^+$ denotes the self-dual part of the curvature $F_A$, and $q(\phi)$ is the trace free part of the endomorphism $\phi\otimes\phi^*$ of $S^+$ and this endomorphism is identified with an imaginary-valued self-dual $2$-form via the Clifford multiplication.

The action of the gauge transformation group $\G=\Map(X;\U(1))$ is given by $u(A,\phi)=(A-2u^{-1}du, u\phi)$ for $u\in\G$. 

Fix $k>4$. 
Let $\CC$ be the $L^2_k$-completion of $\Omega^1(X)\oplus\Gamma(S^+)$, and $\UU$ be the $L^2_{k-1}$-completion of $\Gamma(S^-)\oplus i\Omega^+(X)\oplus\Omega^0(X)/\R$, where $\Omega^k(X)$ is the space of $k$-forms, $\Omega^+(X)$ is the space of self-dual $2$-forms, and $\R$ is the space of constant functions on $X$.
Fix a $\Z_2$-invariant connection $A_0$ on $L$, and let us define the monopole map $\mu\colon \CC\to\UU$ by
$$
\mu (a,\phi) = (D_{A_0+ia}\phi, F^+_{A_0+ia}-q(\phi),d^*a).
$$
When the lift of the $\Z_2$-action to the $\Spinc$-structure $c$ is given, $\mu$ is $\Z_2\times S^1$-equivariant. 

Now we will review on finite dimensional approximations.
Decompose the monopole map $\mu$ into the sum of the linear part $\D$ and the quadratic part $\QQ$, i.e., $\mu = \D + \QQ$, where $\D\colon\CC\to\UU$ is given by 
$$
\D(a,\phi)=(D_{A_0}\phi,d^+a,d^*a), 
$$
and $\QQ$ is the rest.
\begin{Theorem}[\cite{BF}]\label{thm:finite-approx}
There exists a finite dimensional linear subspace $W_f\subset \UU$ which has the following properties{\rom :}
\begin{enumerate}
\item $W_f$ and images of $\D$ spans $\UU$, i.e., $\UU = W_f + \im\D$.
\item For each finite dimensional linear subspace $W\subset \UU$ which contains $ W_f$, let $V= \D^{-1}(W)$.
Then a $\Z_2\times S^1$-equivariant map $f_W\colon S^V\to S^W$ between one-point compactifications of $V$ and $W$ is defined from $\mu$.
{\rom (}The base points of $S^V$ and $S^W$ at infinity are denoted by $*$. {\rom )}
\item If $W^{\prime}=U\oplus W$, where $U$ and $W$ are finite dimensional linear subspaces, and $W\supset W_f$, then $f_{W^{\prime}}$ and $\id_U\wedge f_W$ are $\Z_2\times S^1$-equivariant homotopic as pointed maps{\rom :}
$$
S^{V^\prime}\cong S^{U\oplus V}\to S^{W^\prime}\cong S^{U\oplus W}.
$$ 
\end{enumerate}
\end{Theorem}
\begin{Remark}
More concretely, the map $f_W$ is defined as follows. 
It is known that the monopole map $\mu$ is proper, and extends to the map $\mu^+$ between $S^{\CC}$ and $S^{\UU}$.
For each finite dimensional linear subspace $W$ of $\UU$, let $S(W^\perp)$ be the unit sphere in the orthogonal complement $W^\perp$ of $W$.
In \cite{BF}, a retraction $\rho_W\colon S^{\UU}\setminus S(W^\perp)\to S^W$ is given, and it is proved that $\im \mu^+\cap S(W^\perp)=\emptyset$ when $W$ contains $W_f$.  
Then $f_W = (\rho_W\circ\mu^+|_{S^V})\colon S^V\to S^W$.
\end{Remark}

\thmref{thm:finite-approx} enables us to specify a well-defined class in the following equivariant stable cohomotopy group:
$$
[f_W]\in \{\ind_{\Z_2} D, H^+\}^{\Z_2\times S^1} = \underset{U\subset {\UU}}{\colim}[S^U\wedge S^{V-W}, S^U\wedge S^0]^{\Z_2\times S^1},
$$
where $S^{V-W}$ is the desuspension of $S^V$ by $S^W$ which can be identified with the space of maps from $S^W$ to $S^V$.
Note that the space $\UU$ is the $\Z_2\times S^1$-universe which consists of the following representations: $\C_+$, $\C_-$, $\R_+$, $\R_-$, where $S^1(\subset \C)$ acts on $\C_+$, $\C_-$ by multiplication, and on $\R_+$, $\R_-$ trivially, and $\Z_2$ acts on $\C_+$, $\R_+$ trivially, and on $\C_-$, $\R_-$ by multiplication of $\pm 1$. 
The class $[f_W]$ is the {\it $\Z_2$-equivariant Bauer-Furuta invariant} of $(X,c)$, and we denote it by $\fb^{\Z_2}(c)$.

In the case without any extra group action on the $\Spinc$-structure, a similar construction determines an element $[f_W]$ in $\{\ind D, H^+\}^{S^1}$, and this is the original Bauer-Furuta invariant, denoted by $\fb (c)$.
Note that the original Bauer-Furuta invariant $\fb(c)$ is obtained from $\fb^{\Z_2}(c)$ via the map forgetting the $\Z_2$-action $\phi\colon  \{\ind_{\Z_2} D, H^+\}^{\Z_2\times S^1}\to \{\ind D, H^+\}^{S^1}$, by $\fb(c)=\phi(\fb^{\Z_2}(c))$.

In general, it is difficult to calculate stable cohomotopy groups. 
However, in several cases of $S^1$, stable cohomotopy groups  $\{\ind D, H^+\}^{S^1}$ are calculated in \cite{BF}.
In fact, when the virtual dimension $d(c)$ of the moduli of $c$ is $1$, we know
\begin{equation}\label{eq:d1}
\{\ind D, H^+\}^{S^1} \cong \left\{
\begin{aligned}
 \Z_2, &\text{ when $\ind D$ is even,}\\
 \{0\}, &\text{ otherwise.}
\end{aligned}
\right.
\end{equation}

In the $\Z_2\times S^1$-case, $\{\ind_{\Z_2} D, H^+\}^{\Z_2\times S^1}$ is attained by suspension by finite copies of $\C_+$, $\C_-$ and $\R_+$:
\begin{Proposition}\label{prop:suspension}
Suppose $d(c)=1$ and $b_+\geq 2$.
Let $V_0= k_+\C_+\oplus k_-\C_-$ and $W_0 = H^+$. 
{\rom (}If one of $k_\pm$ is negative, say $k_+$, then add $-k_+\C_+$ to both of  $V_0$ and $W_0$ to obtain actual representations.{\rom )}
There is a representation $V^\prime\subset \UU$ of the form 
$$
V^\prime = a_+\C_+\oplus  a_-\C_-\oplus b\R_+,
$$
such that 
$$
[S^{V^\prime}\wedge S^{V_0}, S^{V^\prime}\wedge S^{W_0}]^{\Z_2\times S^1}\cong \{\ind_{\Z_2} D, H^+\}^{\Z_2\times S^1}.
$$
\end{Proposition}

For the proof, we invoke the equivariant Freudenthal suspension theorem.
\begin{Theorem}[\cite{May}, Chapter IX, Theorem 1.4]\label{thm:Freudenthal}
Let $G$ be a compact Lie group, $U$ a representation, $Y$ a $G$-space, and $X$ a $G$-CW complex {\rom (}see  \secref{sec:eqob}{\rom)}. 
For each subgroup $H\subset G$, let $c^H(Y)$ be the connectivity of $Y^H$. 
If the following hold,
\begin{enumerate}
\item $\dim X^H\leq 2c^H(Y)$ for all subgroups $H$ with $U^H\neq 0$, 
\item $\dim X^H\leq c^K(Y)-1$ for all pairs of subgroups $K\subset H$ with $U^H\neq U^K$,
\end{enumerate}
then the suspension map
$$
S^U\colon [X,Y]^G\to [S^U\wedge X,S^U\wedge Y]^G,
$$
is bijective.
\end{Theorem}
\proof[Proof of \propref{prop:suspension}]
The proof is divided into two steps. 

(1)For each subgroup $H\subset G=\Z_2\times S^1$ such that there is a representation $V\subset \UU$ with $V^H\neq 0$, by adding copies of $V$ to $V^\prime$ if necessary, the following is satisfied:
$$
\dim {V^\prime}^H+\dim V_0^H\leq 2(\dim {V^\prime}^H + \dim W_0^H -1).
$$

(2)For each pair of subgroups $K\subset H \subset G$ such that there is a representation $V\subset \UU$ with $V^H\neq V^K$, by adding copies of $V$ to $V^\prime$ if necessary, the following is satisfied:
$$
\dim {V^\prime}^H+\dim V_0^H\leq (\dim {V^\prime}^K + \dim W_0^K -1)-1.
$$

If (1) and (2) are satisfied, then $[S^{V^\prime}\wedge S^{V_0}, S^{V^\prime}\wedge S^{W_0}]^{\Z_2\times S^1}$ is in stable range by \thmref{thm:Freudenthal}.
Therefore it suffices to prove that we can take a direct sum of finite copies of $\C_\pm$ or $\R_+$ as $V^\prime$.
First, note that $\UU$ contains only four types of representation, and therefore the number of orbit types in them is finite.
For (1), it suffices to add a direct sum of finite copies of $\R_+$ to $V^\prime$.
For (2), we can prove that it is not necessary to add $\R_-$ to $V^\prime$ under the assumption $b_+\geq 2$, and it suffices to add a finite sum of $\C_\pm$ or $\R_+$ to $V^\prime$.
\endproof

Take a large $V^\prime$ as in \propref{prop:suspension}, and put $V = V^\prime\oplus V_0$ and $W = V^\prime\oplus W_0$.
We use the following notation.
Let $V_\R^\pm$ be the $\R_\pm$-component of $V$, and $V_\C^\pm$ the $\C_\pm$-component of $V$, and similarly $W_\R^\pm$ and $W_\C^\pm$.
Let $V_\C= V_\C^+\oplus V_\C^-$, and $V_\R$, $W_\C$ and $W_\R$ are also similarly defined.
%
%
\subsection{Bauer-Furuta invariants as obstruction classes}\label{subsec:obstruction}
%
%
In this subsection, we describe (non-equivariant) Bauer-Furuta invariants in terms of ordinary obstruction theory. First, note the following:
\begin{Proposition}[\cite{BF}, Proposition 3.4]
If $\ind D >0$, then $\{ S^V,S^W\}^{S^1}\cong \{S^V/S^1,S^W\}$.
\end{Proposition}
Note that, if $d(c)=1$, then $\ind D >0$.
Thus, by the Freudenthal theorem, $\{ S^V, S^W \}^{S^1}$ is isomorphic to non-equivariant (unstable) cohomotopy group $[S^V/S^1, S^W]$ for sufficient large $V$, $W$.
Then we can analyse $[S^V/S^1, S^W]$ by the standard obstruction theory.

\begin{Proposition}
Suppose $d(c)=1$ and put $n=\dim S^V/S^1$. {\rom (}Note that $\dim S^W= n-1$.{\rom)}
Then $H^r(S^V/S^1,*\,;\pi_r(S^W)) = 0$ when $r\neq n$, and  $H^n(S^V/S^1,*\,;\pi_n(S^W))\cong \Z/2$. 
\end{Proposition}
\proof
It is clear that $H^r(S^V/S^1,*\,;\pi_r(S^W)) = 0$ when $r< n-1$ or $r>n$.
For sufficient large $V$, the codimension of $S^{V_\R}$ in $S^V/S^1$ is also large.
Now note that $(S^V/S^1,*)$ has the same homotopy type as $(S^{V_\R}\wedge S^1\wedge \CP(V_\C), *)$ ({\it cf.} \eqref{eq:decomp}) and $S^W\cong S^{n-1}$. Then we see the rest case when $r= n-1,n$.
\endproof
Then the standard obstruction theory implies the following:
\begin{Theorem}[See \cite{Hu}, Chapter VI]\label{thm:Hu}
There exists a subgroup $J$ of $H^n(S^V/S^1,*\,;\pi_n(S^W))$  such that
$$
[S^V/S^1, S^W] \cong H^n(S^V/S^1,*\,;\pi_n(S^W))/J.
$$
The isomorphism is given by correspondence that  $f\in [S^V/S^1, S^W]$ is mapped to the difference obstruction class $d(f, \underline{0})$, where $\underline{0}$ is the map which collapse whole $S^V$ to the base point.
\end{Theorem}
\begin{Remark}
The subgroup $J$ in \thmref{thm:Hu} is given as follows.
Let $(K,L)=(S^V/S^1,*\,)$ and $Y=S^W$.
Suppose a map $f\colon K\to Y$ is given. 
Let $K^r$ be the $r$-skeleton of $K$, and put $\bar{K}^r = K^r \cup L$.
Let $\Theta$ be the set of maps $g\colon \bar{K}^{n-1}\to Y$ such that $g|L=f|L$, and put $\theta_0=f|\bar{K}^{n-1}$.
Let $R^n(K,L;f) :=\pi_1(\Theta,\theta_0)$.
Then each element in $R^n(K,L;f)$ can be represented by a homotopy $h_t\colon\bar{K}^{n-1}\to Y$ such that 
$$
h_0=h_1=f|\bar{K}^{n-1},\quad h_t|L=f|L,\quad 0\leq t\leq 1.
$$
Let us define the map $\xi\colon R^n(K,L;f)\to H^n(K,L;\pi_n(S^W))$ by $\xi([h_t]) = d(f,f;h_t)$, which is the obstruction for deforming $f$ to itself via the homotopy $h_t$.
Then $J=\im\xi$.

The meaning of \thmref{thm:Hu} is as follows.
If we take a map, say $\underline{0}$, as origin, then the cohomology $H^n(S^V/S^1,*\,;\pi_n(S^W))$ itself is identified with the set of maps modulo homotopies which are fixed on $(n-2)$-skeleton. 
To obtain the required homotopy set given via homotopies which are not necessarily fixed on $(n-2)$-skeleton, we need to know when two maps whose difference is represented by an element of the cohomology $H^n(S^V/S^1,*\,;\pi_n(S^W))$ are in fact homotopic.
The subgroup $J$ measures this.
Note that \eqref{eq:d1} can be proved by calculating $J$ directly. 
\end{Remark}

Thus, when the map $f_W$ is a finite dimensional approximation of the monopole map, the Bauer-Furuta invariant $\fb(c)$ can be written as the difference obstruction class  $d(f_W,\underbar{0})$. 
%
%
%
%
\section{Equivariant obstruction theory}\label{sec:eqob}
%
%
In some special situations, equivariant Bauer-Furuta invariants can be written as {\it equivariant} obstruction classes. 
In this section, we give a brief review on the equivariant obstruction theory. 
For more details, see \cite{May, Bredon, M2, Dieck}.
%
%
\subsection{Bredon cohomology}\label{subsec:bredon}
%
%
First, we review on Bredon cohomology which is the foundation of the equivariant obstruction theory. 
Let $G$ be a compact Lie group.
We recall the notion of $G$-CW complex. 
\begin{Definition}
A $G$-CW complex $K$ is the union of sub $G$-spaces $K^n$ such that 
\begin{itemize}
\item $K^0$ is a disjoint union of orbits $G/H$,
\item $K^{n}$ is obtained from $K^{n-1}$ by attaching $n$-cells $\sigma \cong G/H_{\sigma}\times D^{n}$ via attaching $G$-maps $a_\sigma\colon G/H_{\sigma}\times S^{n-1}\to K^{n-1}$, 
\item $K$ has the colimit topology with respect to $(K^n)$.
\end{itemize}
\end{Definition}
We denote the maximal degree of cells in $K$ by $\celldim K$. 
In general, $\celldim K= \dim K/G$.
We need several notions to define Bredon cohomology.
\begin{Definition}
The {\it category of canonical orbits} of $G$, denoted by $\OO_{G}$, is the category whose objects are cosets $G/H$, where $H$ are closed subgroups, and morphisms are $G$-homotopy classes of $G$-maps $G/H\to G/I$.
\end{Definition}
\begin{Remark}
Let $N(H,I)$ be the subgroup  of $a\in G$ such that $a^{-1}Ha\subset I$.
Then we have an identification,
$
G\text{-$\map$}(G/H, G/I)\cong N(H,I)/I = (G/I)^H.
$
\end{Remark}
Let $Abel$ be the category of Abelian groups.
\begin{Definition}
We call a contravariant functor $\OO_{G}\to Abel$ a ($G$-){\it coefficient system}.
\end{Definition}
\begin{Remark}
For a closed subgroup $H$, let $WH$ be the Weyl group of $H$, and $W_0H$ its identity component. Then, for a coefficient system $M\colon \OO_{G}\to Abel$, $M(G/H)$ is a $\Z[WH/W_0H]$-module.
\end{Remark}
Let us consider the category of whole $G$-coefficient systems, denoted by $\CC_{G}$, whose morphisms are given by natural transformations.
It is known that $\CC_{G}$ is an Abelian category (\cite{Mac}, p.258).
Therefore we can develop homological algebra in it.

Let $K$ be a $G$-CW complex.
We have coefficient systems, 
$$
\underline{C}_n(K) = \underline{H}_n(K^n,K^{n-1};\Z),
$$
which are defined by $\underline{C}_n(K) (G/H) = H_n((K^n)^H, (K^{n-1})^H;\Z)$.
The connecting homomorphisms of homology of the triples $((K^n)^H, (K^{n-1})^H,(K^{n-2})^H)$ induce a homomorphism of $G$-coefficient systems,
$$
\underline{\partial}\colon \underline{C}_n(K)\to \underline{C}_{n-1}(K),
$$
which satisfies $\underline{\partial}\circ \underline{\partial}=0$.

Let $M$ be a $G$-coefficient system, and define the cochain complex $(C_{G}^n(K;M),\delta)$ by
$$
C_{G}^n(K;M) = \Hom_{\CC_{G}}(\underline{C}_n(K) ,M),\quad \delta = (\underline{\partial})^*.
$$
Its cohomology is the {\it Bredon cohomology}, denoted by $H_{G}^n(K;M)$.

The coefficient systems of chain, $\underline{C}_n(K)$, can be identified as follows.
\begin{Proposition}
We have an isomorphism,
$$
\underline{C}_n(K)\cong \bigoplus_{\sigma: \text{ $n$-cell}} \underline{H}_0 (G/H_\sigma;\Z),
$$
where $\underline{H}_0 (G/H_\sigma;\Z)$ is defined by $\underline{H}_0 (G/H_\sigma;\Z)(G/H) = H_0((G/H_\sigma)^H;\Z).$
\end{Proposition}
\proof
This follows from
$$
H_n((K^n)^H,(K^{n-1})^H)\cong \bigoplus_{\sigma\text{: $n$-cell}} H_n((G/H_\sigma)^H \times (D^n,S^{n-1}))\cong \bigoplus_{\sigma\text{: $n$-cell}} H_0((G/H_\sigma)^H).
$$
\endproof
\begin{Remark}
Note  that
$$
H_0((G/H_\sigma)^H;\Z)\cong F\pi_0((G/H_\sigma)^H)\cong F[G/H,G/H_\sigma ]^{G},
$$
\end{Remark}
where, for a set $S$, $FS$ denotes the free Abelian group generated by $S$.

The coefficient system $\underline{H}_0(G/H)$ has a special property in the next.
\begin{Proposition}\label{prop:homcg}
Let $M$ be a coefficient system. There is an isomorphism,
$$
\Hom_{\CC_{G}}(\underline{H}_0(G/H),M)\cong M(G/H),
$$
via $\phi\mapsto \phi(1_{G/H})$.
\end{Proposition}
\proof
The inverse is given as follows.
Suppose that $m_\phi\in M(G/H)$ to be $\phi(1_{G/H})$ is given.
Recall that $\underline{H}_0 (G/H;\Z)(G/I)\cong F[G/I,G/H ]^{G}$. 
For a $G$-map $f\colon G/I\to G/H$, $M(f)m_\phi$ is in $M(G/I)$.
Then, the correspondence $m_\phi\mapsto M(f)m_\phi$ determines a morphism $\phi$ in $\CC_{G}$.
\endproof
\begin{Corollary}
The cochain group $C_{G}^n(K;M)$ can be identified with the space of functions such that an $n$-cell $\sigma$ is mapped to an element of $M(G/H_\sigma)$.
\end{Corollary}
\begin{Remark}
To define $H_{G}^n(K;M)$, the coefficient system $M$ need not to be defined on the whole $\OO_G$, but on a subcategory $\OO^K_G$ given below. 
Let $\Iso(K)$ be the set of isotropy groups in $K$. 
Objects of $\OO^K_G$ are cosets $G/H$ where $H\in\Iso(K)$. Morphisms are $G$-homotopy classes of $G$-map between them.
\end{Remark}

For a $G$-CW pair $(K,L)$, the relative Bredon cohomology $H_{G}^*(K,L;M)$ is defined from
$$
\underline{C}_n (K,L)= \underline{H}_n (\bar{K}^n,\bar{K}^{n-1};\Z),
$$
where $\bar{K}^k = K\cup L$.
Also there is an identification $\underline{C}_n (K,L)=\bigoplus_{\sigma}\underline{H}_0(G/H_\sigma)$, where $\sigma$ runs over $n$-cells in $K\setminus L$.

This definition works well due to the following.
\begin{Proposition}\label{prop:proj}
For a $G$-CW pair $(K,L)$, $\underline{C}_n (K)$, $\underline{C}_n (L)$ and $\underline{C}_n (K,L)$ are projective in $\CC_{G}$.  
\end{Proposition}
This is proved from the fact that the coefficient system $\underline{H}_0(G/H)$ is projective, which can be easily seen from \propref{prop:homcg}. 

For a $G$-CW pair $(K,L)$, we have an exact sequence,
$$
0\to \underline{C}_n (L)\to \underline{C}_n (K)\to \underline{C}_n (K,L)\to 0.
$$
When we apply the functor $\Hom_{\CC_{G}}(\,\cdot\,, M)$ to this, the resulting sequence is also exact by \propref{prop:proj}.
Thus we obtain long exact sequences of Bredon cohomology.

In general, it is not easy to calculate Bredon cohomology.
However, the following simple case is rather easy. 
\begin{Proposition}[{\it Cf.} \cite{Dieck}, Chapter \Romnum{2}, 3]\label{prop:free}
Let $G$ be a compact Lie group, and $G_0$ its identity component.
Let $M$ be a $G$-coefficient system, and put $M_0 = M(G/\{e\})$. 
{\rom (}Note that $M_0$ is a $\Z[G/G_0]$-module.{\rom )}
Suppose that $(K,L)$ is a $G$-CW pair such that $G$ acts on $K\setminus L$ freely. 
Let us consider the complex{\rom :}
\begin{align*}
\to
\Hom_{\Z[G/G_0]}(H_{n-1}(\bar{K}^{n-1}/G_0, \bar{K}^{n-2}/G_0),M_0)\overset{\delta_{n-1}}{\to} \Hom_{\Z[G/G_0]}(H_{n}(\bar{K}^{n}/G_0, \bar{K}^{n-1}/G_0),M_0)\\ 
\overset{\delta_n}{\to}\Hom_{\Z[G/G_0]}(H_{n+1}(\bar{K}^{n+1}/G_0, \bar{K}^{n}/G_0),M_0)\to,
\end{align*}
where $\delta_r$ are determined from connecting homomorphisms of homology of triples 
$$
(\bar{K}^{r}/G_0, \bar{K}^{r-1}/G_0,\bar{K}^{r-2}/G_0).
$$
Then the $n$-th Bredon cohomology $H_{G}^n(K,L;M)$ is isomorphic to $\ker\delta_{n}/\im\delta_{n-1}$.
\end{Proposition}
\proof
By the assumption, we have $\underline{C}_n(K,L)\cong \bigoplus_{\sigma}\underline{H}_0(G/\{e\})$. 
Then note that $H_0(G)\cong H_0(G/G_0)$ and \propref{prop:homcg}.
\endproof
\begin{Remark}
If $(K,L)$ is a $G$-CW complex, then $(K/G_0,L/G_0)$ has an ordinary CW complex structure with a cellular $G/G_0$-action induced from the $G$-CW structure on $(K,L)$. 
In the case of \propref{prop:free}, we can calculate the Bredon cohomology of $(K,L)$ from this CW complex structure of $(K/G_0,L/G_0)$ with the $G/G_0$-action. 
\end{Remark}
%
%
\subsection{Equivariant obstruction theory}\label{subsec:eqob}
%
%
We can develop equivariant obstruction theory on Bredon cohomology.
Let $(K,L)$ be a $G$-CW pair. 
Let $\Iso(K,L)$ be the set of isotropy groups of $K\setminus L$.
\begin{Definition}
Fix $n\geq 1$.
Let $Y$ be a $G$-space such that, for each closed subgroup $H\in\Iso(K,L)$, $Y^H$ is non-empty, path-connected and $n$-simple.
Then the coefficient system $\underline{\pi}_n(Y)$ is defined by $\underline{\pi}_n(Y)(G/H)=\pi_n(Y^H)$.
\end{Definition}
Suppose that a $G$-map $f\colon \bar{K}^n = K^n\cup L\to Y$ is given. 
We ask when $f$ can be extended $G$-equivariantly over $\bar{K}^{n+1}$.
For each $(n+1)$-cell $\sigma$, the composite map of the attaching map $a_\sigma\colon G/H_\sigma\times S^n\to \bar{K}^n$ with $f$ defines an element of $\pi_n(Y^{H_\sigma})$ since $f\circ a_{\sigma}$ maps $eH_\sigma\times S^n$  into $Y^{H_\sigma}$.
This correspondence gives a well-defined cocycle,
$$
c_f\in C^{n+1}_{G}(K,L;\underline{\pi}_n(Y)),
$$
and $f$ extends to $\bar{K}^{n+1}$ if and only if $c_f=0$.

Similarly, for $G$-maps $f_0,f_1\colon \bar{K}^n\to Y$ and a $G$-homotopy $h_t$ on $\bar{K}^{n-1}$ rel $L$ such that $h_0=f_0|\bar{K}^{n-1}$ and $h_1=f_1|\bar{K}^{n-1}$,  we can define the difference cochain, 
$$
d(f_0,f_1;h_t)\in C^{n}_{G}(K,L;\underline{\pi}_n(Y)),
$$
such that $\delta d(f_0,f_1;h_t) = c_{f_1} - c_{f_0}$.
Then $h_t$ extends over $\bar{K}^{n}$ if and only if $d(f_0,f_1;h_t)=0$.
Furthermore, for given $f_0\colon \bar{K}^n\to Y$ and $\beta\in  C^{n}_{G}(K,L;\underline{\pi}_n(Y))$, there exists $f_1\colon \bar{K}^n\to Y$ which coincides with $f_0$ on $\bar{K}^{n-1}$ and satisfies $\beta = d(f_0,f_1)$, where the homotopy is assumed constant.  

With these understood, a standard argument proves the following.
%
%
%
%
\begin{Theorem}\label{thm:eqHu}
Let $n=\celldim (K\setminus L)$, and denote the set of $G$-homotopy classes rel $L$ of $G$-maps by $[K,Y]_{\rel L}^{G}$.
If $H_{G}^k(K,L;\underline{\pi}_k(Y))=0$ for $k\neq n$, then there exists a subgroup $J^\prime\subset H_{G}^n(K,L;\underline{\pi}_n(Y))$, and 
$$
[K,Y]_{\rel L}^{G}\cong H_{G}^n(K,L;\underline{\pi}_n(Y))/J^\prime.
$$
\end{Theorem}
As in \subsecref{subsec:eqbf}, let $f_W\colon S^V\to S^W$ be a finite dimensional approximation for sufficiently large $V$ and $W$.
If $(K,L)= (S^V,*\,)$ and $Y=S^W$, and the conditions in \thmref{thm:eqHu} are satisfied, then the equivariant Bauer-Furuta invariant $\fb^{\Z_2}(c)$ can be identified with the difference obstruction class $d(f_W,\underline{0})\in H_{\Z_2\times S^1}^n(K,L;\underline{\pi}_n(Y))/J^\prime$.
%
%
%
%
\section{Proof of \thmref{thm:main}}\label{sec:main}
%
%
The purpose of this section is to prove \thmref{thm:main}. 
The proof of \thmref{thm:main} is divided into two steps.
The first step proves Bredon cohomology groups vanish in all degrees except the top . 
Then, by \thmref{thm:eqHu}, $\{\ind_{\Z_2} D, H^+\}^{\Z_2\times S^1}$ is written in terms of the top-degree Bredon cohomology.
The second step compares the top-degree Bredon cohomology with the top-degree ordinary cohomology via the map forgetting the group action, and observes that the forgetting map is the $0$-map.
Then, it follows that the vanishing of the Bauer-Furuta invariant, since this is an image of the forgetting $0$-map. 
Our argument is analogous to that in \cite{Dieck},~Chapter~\Romnum{2},~4 as well as \cite{B2}.
%
%
\subsection{Vanishing of Bredon cohomology in low degrees}\label{subsec:low}
%
%
Recall that we took sufficiently large $V$ and $W$ in \subsecref{subsec:eqbf}.
We assume a $\Z_2\times S^1$-CW complex structure on $(S^V,*\,)$ fixed.
(Later, we will give a $\Z_2\times S^1$-CW complex structure on $(S^V,*\,)$ concretely.)

Note that 
$$
\Iso(S^V,*\,) =\left\{\{e\},\, \Z_2,\,\tilde{\Z}_2,\, \{1\}\times S^1,\, \Z_2\times S^1\right\}.
$$
where, by assuming $\Z_2=\{\pm 1\}$ and $S^1\subset \C$, $\Z_2$ is the subgroup generated by $(-1,1)$, and $\tilde{\Z}_2$ generated by $(-1,-1)$.
\begin{Lemma}
For each $H\in\Iso(S^V,*\,)$, $(S^W)^H$ is a $k$-sphere, where $k\geq 1$.
\end{Lemma}
\proof
This follows from the fact that  $W^H$ is a linear subspace of $W$, $W^H\supset W^{\Z_2\times S^1}$ and $\dim W^{\Z_2\times S^1}\geq b_+^{\Z_2}\geq 1$.
\endproof
Therefore, $\underline{\pi}_n(S^W)$ is well-defined.

Now, we prove the vanishing of Bredon cohomology in low degrees.
\begin{Lemma}
Let $n = \celldim S^V$. 
If $2k_{\pm}< 1 + b_+^{\Z_2}$, then 
$C_{\Z_2\times S^1}^k (S^V,*\,;\underline{\pi}_k(S^W))=0$ for $k\leq n-2$. 
\end{Lemma}
\proof
Let $\sigma$ be a $k$-cell, and $\varphi$ be a $k$-cochain.
Put $G = \Z_2\times S^1$.
If $H_\sigma = \{e\}$, then $\varphi(\sigma)=0$ since $\underline{\pi}_k(S^W)(G/\{e\}) \cong \pi_k(S^{n-1}) =0$.

If $H_\sigma = \{1\}\times S^1$, then $k\leq\dim V_\R$. Since $(S^W)^{S^1} = S^{W_\R}$, $W_\R\cong V_\R\oplus H^+$ and $b_+ >0$, we have $\underline{\pi}_k(S^W)(G/S^1) =0$. 
Similarly, in the case when  $H_\sigma = \Z_2\times S^1$, we can show $\underline{\pi}_k(S^W)(G/G) =0$ by $b_+^{\Z_2} >0$.

Now suppose that $H_\sigma = \Z_2$.
Note that the condition $2k_+<1 + b_+^{\Z_2}$ is equivalent to the condition 
$$
\dim (S^V)^{\Z_2} < 1 + \dim (S^W)^{\Z_2}.
$$
On the other hand, $k\leq \celldim (S^V)^{\Z_2} = \dim (S^V)^{\Z_2} -1$.
Therefore $\underline{\pi}_k(S^W)(G/\Z_2) =0$.

The case when $H_\sigma = \tilde{\Z}_2$ is similar.
\endproof
%
%
\subsection{Calculation in high degrees}\label{subsec:high}
%
%
%
Let $(S^V)_s$ be the singular part of $S^V$, that is,
$$
(S^V)_s=\bigcup_{\{e\}\neq H\in\Iso(S^V,*\,)}(S^V)^H.
$$
\begin{Lemma}\label{lem:codim}
If necessary, by adding several copies of $U = \C_+\oplus\C_-$ to $V$ and $W$, the following holds.
\begin{equation}\label{eq:codim}
\dim (S^V)_s\leq \dim S^V - 4.
\end{equation}
\end{Lemma}
\proof
This follows from the fact that $U^{\Z_2} = \C_+\oplus \{0\}$, $U^{\tilde{\Z}_2} = \{0\}\oplus \C_-$ and these subspaces have complex codimension $1$ in $U=U^{\{e\}}$. 
(Note that $U^{S^1}= U^{\Z_2\times S^1} = \{0\}$.)
\endproof
Let $n=\celldim S^V$. Then, \lemref{lem:codim} and the long exact sequence imply that
$$
H^k_{\Z_2\times S^1}(S^V,*\,;M)\cong H^k_{\Z_2\times S^1}(S^V,(S^V)_s;M),
$$
for $k=n,n-1$.
Therefore, by \propref{prop:free}, $H^k_{\Z_2\times S^1}(S^V,*\,;M)$ for $k=n,n-1$ can be calculated from $S^V/S^1$.

Note that $S^V/S^1$ decomposes as
\begin{equation}\label{eq:decomp}
S^V/S^1 = * \cup V_\R\times(\R_{>0}\times\CP(V_\C)\cup\{0\}).
\end{equation}
Recall that we may assume $V_\R$ contains no $\R_-$ by \propref{prop:suspension}.
Therefore, by giving a $\Z_2$-CW complex structure on $\CP(V_\C)$, we can fix that on $S^V/S^1$ via the decomposition \eqref{eq:decomp}.

Now, we give a $\Z_2$-CW complex structure on $\CP(V_\C)$.
For $a\C_+\oplus b\C_-$, we denote  $\CP(a\C_+\oplus b\C_-)$ by $P^{a,b}$.
We have a filtration as $
P^{a,b}\supset P^{a-1,b}\supset P^{a-2,b}$ if $a\geq 2$. Note that 
$$
P^{a,b}\setminus P^{a-1,b} =\{[z_1,\ldots,z_a,w_1,\ldots,w_b]\,|\, z_a\neq 0\}\cong (a-1)\C_+\oplus b\C_-.
$$
Certainly, this is an open disk.
However, $\Z_2$ still acts on it non-trivially. 
Hence, we need to take a subdivision. It suffices to subdivide $\C_-^b$.

The space $\C_-^b$ is identified with $\R^{2b}=\{(x_1,y_1,\ldots,x_b,y_b)\}$ by $w_i=x_i+\sqrt{-1}y_i$. 
Then, a $\Z_2$-equivariant subdivision is given as follows:
$$
\{x_1>0\}\cup \{x_1<0\}\cup\{x_1=0,y_1>0\}\cup\{x_1=0,y_1<0\}\cup\cdots.
$$
Let us denote the open disk with the above subdivision by $D^{a-1,b}$. 
We may assume $a,b\geq 2$, and have the decomposition:
$$
P^{a,b} = D^{a-1,b}\cup D^{a-2,b}\cup P^{a-2,b}.
$$
According to the above cell structure, we can write down its (ordinary) chain complex $\Z_2$-equivariantly.
$$
\begin{CD}
0@>>> C_n @>\partial_n>> C_{n-1}@>\partial_{n-1}>>  C_{n-2}@>\partial_{n-2}>>  C_{n-3}@>>> \cdots\\
 @.\Vert@.\Vert@.\Vert@.\Vert@.\\
 @.\Z[\Z_2]@.\Z[\Z_2]@.\Z[\Z_2]\oplus\Z[\Z_2]@.\Z[\Z_2]\oplus\Z[\Z_2]@.
\end{CD}
$$
where $\partial_n$ is  given by the matrix over $\Z$ as follows,
\begin{equation}\label{eq:boundary}
\partial_n = 
\begin{pmatrix}
 1 & -1 \\
-1 &  1
\end{pmatrix}.
\end{equation}
\begin{Proposition}
If $b_+-b_+^{\Z_2}$ is odd, then $H^{n-1}_{\Z_2\times S^1}(S^V, (S^V)_s;\underline{\pi}_{n-1}(S^W))=0$.
\end{Proposition}
\proof
When $b_+-b_+^{\Z_2}$ is odd, the $\Z_2$-action on $S^W$ reverses the orientation of $S^W$. 
Therefore the $\Z_2$-action on $\pi_{n-1}(S^W)\cong\Z$ is non-trivial.
With this understood, we use \propref{prop:free}.
In fact, by \eqref{eq:boundary}, we have $\ker\delta_{n-1}=0$ on the $(n-1)$-th cochain group. 
\endproof
\begin{Remark}
If $b_+-b_+^{\Z_2}$ is even, then $H^{n-1}_{\Z_2\times S^1}(S^V, (S^V)_s;\underline{\pi}_{n-1}(S^W))\cong \Z_2$.
\end{Remark}
%
%
\subsection{Completion of the proof}\label{subsec:proof}
%
%
So far, under the assumptions of \thmref{thm:main}, we have the following:
\begin{align*}
H^k(S^V/S^1,*\,;\pi_k(S^W))=0\text{ for }k\neq n,&\quad\{S^V,S^W\}^{S^1}\cong H^n(S^V/S^1,*\,;\pi_n(S^W))/J,\\
H^{k}_{\Z_2\times S^1}(S^V,*\,;\underline{\pi}_{k}(S^W))=0\text{ for }k\neq n,&\quad \{S^V,S^W\}^{\Z_2\times S^1}\cong H^{n}_{\Z_2\times S^1}(S^V,*\,;\underline{\pi}_{n}(S^W))/J^\prime.
\end{align*}

Let us compare $H^n(S^V/S^1,*\,;\pi_n(S^W))$ and $H^{n}_{\Z_2\times S^1}(S^V,*\,;\underline{\pi}_{n}(S^W))$.
Note that the following commutative diagram:
$$
\begin{CD}
0@<<< H^n(S^V/S^1,*\,;\pi_n(S^W))@<<<\Hom_\Z(C_n,\Z_2)@<\delta<<\Hom_\Z(C_{n-1},\Z_2)\\
@. @A{\phi_*}AA @A{\phi}AA @A{\phi}AA \\
0@<<<H^{n}_{\Z_2\times S^1}(S^V,\,*\,;\underline{\pi}_{n}(S^W))@<<<\Hom_{\Z[\Z_2]}(C_n,\Z_2)@<\delta<<\Hom_{\Z[\Z_2]}(C_{n-1},\Z_2),
\end{CD}
$$
where $\phi$ is the map forgetting the $\Z_2$-action. 
Then, a direct calculation shows,
\begin{equation}\label{eq:zero}
\im \phi_* =0.
\end{equation}
We have another commutative diagram,
$$
\begin{CD}
H^n(S^V/S^1,*\,;\pi_n(S^W))@>>> H^n(S^V/S^1,*\,;\pi_n(S^W))/J\cong \{S^V,S^W\}^{S^1}\\
@A{\phi_*}AA @A{\bar{\phi}_*}AA\\
H^{n}_{\Z_2\times S^1}(S^V,*\,;\underline{\pi}_{n}(S^W))@>>> H^{n}_{\Z_2\times S^1}(S^V,*\,;\underline{\pi}_{n}(S^W))/J^\prime\cong \{S^V,S^W\}^{\Z_2\times S^1}.
\end{CD}
$$
Recall that $\fb (c)=\bar{\phi}_*(\fb^{\Z_2}(c))$ and  $\phi_*$ is the $0$-map.
Then, the diagram implies $\fb (c)=0$.
Thus, \thmref{thm:main} is established.
%
%
\subsection{Remarks}\label{subsec:remarks}
%
%
Before ending this section, we give several remarks.
\begin{Remark}
We can give an alternative proof of the mod $p$ vanishing theorem of ordinary Seiberg-Witten invariants in \cite{Fang,Np} in the case when $b_1=0$ by the same method as in this section. 
\end{Remark}
\begin{Remark}
In the case when $d(c)=1$, suppose we are given a $\Z_p$-action of odd prime $p$ instead of a $\Z_2$-action. 
Then some part of arguments in this section also works, however some part does not. 
For example, under appropriate assumptions on $\Z_p$-index of the Dirac operator, we can prove the vanishing of Bredon cohomology groups in low degrees.
However, even if the $(n-1)$-th Bredon cohomology vanishes, \eqref{eq:zero} does not hold.
Therefore, we can not expect such a vanishing result of Bauer-Furuta invariants.
\end{Remark}
\begin{Remark}
Also, in the cases when $d(c)\geq 2$, it would not be easy to prove such a vanishing result, since $(n-2)$-th ordinary cohomology does not vanish.
\end{Remark}
\begin{Remark}\label{rem:odd}
In \thmref{thm:main}, we supposed that the $\Z_2$-action on $X$ lifts to a $\Z_2$-action on $c$. 
Now, let us suppose that the $\Z_2$-action on $X$ only preserves the $\Spinc$-structure $c$, i.e., $\iota^*c\cong c$ for the involution $\iota$ generating the $\Z_2$-action.
Then, $\iota$ also preserves the determinant line bundle $L$. 
Therefore the $\Z_2$-action lifts to $L$ by the theorem of Hattori and Yoshida~\cite{HY}. 
Let $P_{\SO}$ be the frame bundle of $X$ and $P_{\U(1)}$ the $\U(1)$-bundle for $L$. 
Recall that the $\Spinc$-structure $c$ is given by a double covering $P_{\Spinc(4)}\to P_{\SO}\times_X P_{\U(1)}$. 
Since the $\Z_2$-action on $P_{\SO}\times_X P_{\U(1)}$ is given, two cases may occur on the lifting of the $\Z_2$-action to $c$.
The first case is that the $\Z_2$-action on $X$ lifts to a $\Z_2$-action on $c$ as above. 
We call such a case that the $\Z_2$-action is {\it even type} with respect to $c$.
The second case is that the $\Z_2$-action on $X$ does not lift as a $\Z_2$-action on $c$, however it is covered by a $\Z_4$-action on $c$. 
We call such a case that the $\Z_2$-action is {\it odd type} with respect to $c$.

In the odd case, the $\Z_4$-index of the Dirac operator is written as $\ind_{\Z_4} D_{A_0}= k_1\C_1 +k_3\C_3$, where $\C_j$ is the complex $1$-dimensional representation of weight $j$. 
(This is because the $\Z_4$-lift of the generator of $\Z_2$ acts on spinors as multiplication by $\pm\sqrt{-1}$.)

In the odd case, we can also prove the following result similar to \thmref{thm:main}. 
\end{Remark}
\begin{Theorem}\label{thm:odd}
Suppose the following conditions are satisfied{\rom :}
\begin{enumerate}
\item $b_1=0$, $b_+\geq 2$ and $b^{\Z_2}_+\geq 1$,
\item $b_+-b_+^{\Z_2}$ is odd,
\item $d(c)=2(k_1 + k_3) - (1 + b_+)=1$,
\item $2k_{j}< 1+ b_+^{\Z_2}$, for $j=1,3$.
\end{enumerate}
Then the Bauer-Furuta invariant of $c$ is zero.
\end{Theorem}
The proof of \thmref{thm:odd} is analogous to that of \thmref{thm:main}.
In the odd case, the finite dimensional approximation $f_W\colon V\to W$ is $\Z_4\times S^1$-equivariant, and $V$ and $W$ are direct sums of finite copies of $\C_1$, $\C_3$, $\R_\pm$. 
Note that every point in $V$, $W$ have an isotropy which contains a subgroup $\Z_2=\langle(g^2,-1)\rangle$, where $g$ is a generator of $\Z_4$.
Let $G=\Z_4\times S^1/\langle(g^2,-1)\rangle\cong \Z_2\times S^1$. 
Then $f_W$ is $G$-equivariant.
Now, an argument  similar as in previous subsections proves \thmref{thm:odd}.
%
%
%
%
\section{Applications}\label{sec:application}
%
%
In this section, we prove \thmref{thm:nonsmoothable} as an application of \thmref{thm:main}.
Below, the Euler number of a manifold $Y$ is denoted by $\chi(Y)$, and the signature by $\Sign(Y)$.
We assume every manifold is oriented and every $\Z_2$-action is orientation-preserving, unless stated otherwise.
%
%
\subsection{Constraint on smooth actions on $K3\#K3$}\label{subsec:const}
%
%
Rewriting the condition on the $\Z_2$-index of the Dirac operator in \thmref{thm:main} by the $G$-index theorem, we can relate fixed point data to the value of the Bauer-Furuta invariant.

Let $X$ be a simply-connected smooth spin $4$-manifold with a $\Z_2$-action.
When $X$ is simply-connected, the spin structure on $X$ is unique up to equivalence.  
Therefore every involution $\iota\colon X\to X$ preserves the spin structure and also the $\Spinc$-structure $c_0$ which is determined by the spin structure. 
Then, as mentioned in \remref{rem:odd}, $\iota$ lifts to $c_0$ as a $\Z_2$ (even type) or $\Z_4$ (odd type) action. 

By \cite{AB}, the $\Z_2$-action is even type if and only if the fixed point set is discrete.
Now, for an even-type  action, the $G$-spin theorem \cite{AB} claims that
\begin{align*}
\ind_g D &= k_+ - k_- = \frac14 \sum_{p\in X^{\Z_2}} \varepsilon(p),\\
\ind D &= k_+ + k_- = -\frac18 \Sign (X),
\end{align*}
where $\varepsilon\colon X^{\Z_2}\to \{\pm 1\}$ is the sign assignment determined by the lift of the $\Z_2$-action to $c_0$. 
By solving the above equations, we obtain
\begin{equation}\label{eq:G-ind}
2k_{\pm} = -\frac18 \Sign (X)\pm  \frac14 \sum_{p\in X^{\Z_2}} \varepsilon(p).
\end{equation}
Note that both $k_+$ and $k_-$ are even because of the quartanion structure of Dirac index.
Note also that the sum $\sum \varepsilon(p)$ is a multiple of $8$, because $\frac18 \Sign (X)$ is even.

If the  Bauer-Furuta invariant of  $c_0$ is non-trivial, then \thmref{thm:main} implies a constraint on the sign assignment $\varepsilon$. 
Furuta, Kametani and Minami proved that a spin manifold which has the same rational cohomology ring as $K3\#K3$ has non-trivial Bauer-Furuta invariant for every spin structure~\cite{FKM}.
Therefore, we obtain the following.
\begin{Proposition}\label{prop:K3K3}
Let $X$ be a smooth closed simply-connected spin manifold which has the same intersection form as $K3\#K3$.
Suppose an orientation-preserving smooth $\Z_2$-action on $X$ is given, and let $c_0$ be the $\Spinc$-structure which is determined by a spin structure.
If $b_+^{\Z_2}=5$ and fixed points are discrete, then the lift of $\Z_2$-action to $c_0$ as above satisfies, 
\begin{equation}\label{eq:epsilon}
\sum_{p} \varepsilon(p)=\pm  8.
\end{equation}
\end{Proposition}
\proof
Let $m$ be the number of $X^{\Z_2}$. 
First, we prove $m=10$ if $b_+^{\Z_2}=5$. 
By the ordinary Lefschetz formula, we have $\chi(X/\Z_2) = \frac12(\chi(X)+m)$.
On the other hand, the $G$-signature theorem implies that $\Sign(X/\Z_2)=\frac12 \Sign(X)$.
From the equation $1+b_+^{\Z_2}=\frac12(\chi(X/\Z_2)+ \Sign(X/\Z_2))$, we have $m=10$.
Since $\sum\varepsilon(p)$ is a multiple of $8$, we have $\sum\varepsilon(p)= -8$, $0$ or $8$.
If $\sum\varepsilon(p)= 0$, then $2k_{\pm}=4<6=1+b_+^{\Z_2}$. 
Therefore \thmref{thm:main} implies $\fb(c_0)=0$.
However this contradicts with the non-vanishing result for $K3\#K3$ in \cite{FKM}.
\endproof
%
%
\subsection{Atiyah-Bott's criterion for $\varepsilon(p)$}\label{subsec:criterion}
%
%
Before going further, let us recall Atiyah-Bott's criterion for $\varepsilon(p)$. 
Let $X$ be a spin $4$-manifold, and suppose that a smooth involution $\iota\colon X\to X$ with isolated fixed points is given.
When an $\iota$-invariant metric is fixed,  $\iota$ lifts to the frame bundle $F$ as $\iota_*\colon F\to F$.
A spin structure on $X$ is given by a double cover $\varphi\colon\hat{F}\to F$, where $\hat{F}$ is a $\Spin(4)$ bundle.
Suppose that $\iota$ lifts to $\hat{F}$.

For distinct fixed points $P$ and $Q$, we shall compare $\varepsilon(P)$ and $\varepsilon(Q)$.
For this purpose, we take a path $s$ in $F$ starting from a point $y\in F_P$ and ending at $y^\prime \in F_Q$.   
Then the path $-\iota_*s$ has the same starting point and the end point as $s$, where ``$-$'' means the multiplication by $-1$ on each fiber.
Thus, by connecting $s$ with $-\iota_*s$, we obtain a circle $C$ in $F$.

Then, Atiyah-Bott's criterion is given as follows.
\begin{Proposition}[\cite{AB}]\label{prop:criterion}
The preimage $\varphi^{-1}(C)$ has two components if and only if $\varepsilon(P) = \varepsilon(Q)$. 
In other words, the preimage $\varphi^{-1}(C)$ is connected if and only if $\varepsilon(P) = -\varepsilon(Q)$. 
\end{Proposition}
%
%
%
\subsection{Edmonds-Ewing's construction of locally linear $\Z_2$-actions}\label{subsec:EE}
%
%
To construct locally linear $\Z_2$-actions, we invoke the realization theorem by Edmonds and Ewing~\cite{EE}. 
We summarize their result in a very special case for our purpose.
\begin{Theorem}[\cite{EE}]\label{thm:EE}
Suppose that we are given a $\Z_2$-invariant bilinear unimodular even form $\Psi\colon V\times V\to \Z$ which satisfies the following{\rom :}
\begin{enumerate}
\item As a $\Z[\Z_2]$-module, $V\cong T\oplus F$, where $T$ is a trivial $\Z[\Z_2]$-module with $\rank_{\Z}T=n$, and $F$ is a free $\Z[\Z_2]$-module. 
\item For any $v\in V$, $\Psi (v,gv)\equiv 0$ mod $2$, where $g$ is the generator of $\Z_2$.
\item The $G$-signature theorem is satisfied, i.e., $\Sign(g, (V,\Psi))=0$.
\end{enumerate}
Then, there exists a locally linear $\Z_2$-action on a simply-connected $4$-manifold $X$ such that its intersection form is $\Psi$, and the number of fixed points is $n+2$.
\end{Theorem}
\begin{Remark}
Since the form $\Psi$ is assumed even in \thmref{thm:EE}, the homeomorphism type of $X$ is unique by Freedman's theorem.
\end{Remark}
For our application, we need to recall their construction precisely. 
Their construction is an equivariant handle construction.

Let $B_0$ be a unit ball in $\C^2$ on which $\Z_2$ acts by multiplication of $\pm 1$.
Let us take a $\Z_2$-invariant knot $K$ in $S_0=\partial B_0$.
Then a framing of $K$ can be represented  by an equivariant embedding $S^1\times D^2\to S_0$ for some $\Z_2$-action on $S^1\times D^2$. 
In particular, $0$-framing is represented by $f_0\colon S^1\times D^2\to S_0$ such that $f_0(S^1\times\{0\}) = K$, and $f_0(S^1\times\{1\})$ has linking number $0$ with $K$, and the $\Z_2$-action on  $S^1\times D^2\subset \C^2$ is given by $g(z,w)= (-z,-w)$.
An arbitrary $r$-framing of $K$ can be represented by a map $f_r\colon S^1\times D^2\to S_0$ given by $f_r(z,w)=f_0(z,z^rw)$. 
Then $f_r$ is equivariant if the $\Z_2$-action on  $S^1\times D^2$ is given by $g(z,w)= (-z,(-1)^{r-1}w)$.

For a given $K$ and a framing $r$, we can construct a $4$-manifold with a $\Z_2$-action as $W=B_0\cup_{f_r}D^2\times D^2$.

Let $H_1,\ldots,H_n$ be copies of $D^2\times D^2$ on which $\Z_2$ acts by $g(z,w)=(-z,-w)$. 
Note that, if the framing $r$ is even, then we can attach $H_i$ to $B_0$ equivariantly via $f_r$.

Edmonds-Ewing's construction of locally linear actions is divided into three steps.

\smallskip\noindent
{\it \underline{Step 1}}.
Under the assumption of \thmref{thm:EE}, by changing basis of $T$,  we may assume $\Psi|T$ is represented by a matrix $(a_{ij})$  such that $a_{ii}$ is even and $a_{ij}$ is odd whenever $i\neq j$. 
(See \cite{EE}, Section $6$.)
Then we can take a $n$-component link $L_T$ in $S_0 = \partial B_0$ representing the matrix $(a_{ij})$ such that each  component of $L_T$ is $\Z_2$-invariant as follows.
For every two $\Z_2$-invariant knots $K$ and $K^\prime$ in $S_0$, the linking number $lk(K,K^\prime)$ is odd by \cite{EE}, Lemma 5.1. 
Take an arbitrary $n$-component link $L_0$.
By moving various components of $L_0$ across various components of $L_0$ {\it equivariantly}, we can alter the off-diagonal entries of the linking matrix of $L_0$ by arbitrary multiples of $2$. 
Then we can arrange the link to represent the given matrix $(a_{ij})$.
Thus $\Psi |T$ is represented by a framed link $L_T$ in $S_0$.

Now it is not difficult to realize the other part of $\Psi$ by a link, and therefore we obtain a framed link $L$ in $S_0$ which realizes the given $\Z_2$-invariant form $\Psi$.

\smallskip\noindent
{\it \underline{Step 2}.}
Since the diagonal entries of $(a_{ij})$ are assumed even, we can attach $H_1,\ldots, H_n$ and free $2$-handles to $B_0$ along $L$ equivariantly. 
Thus we obtain a $4$-manifold $X_0$ with a smooth $\Z_2$-action,
$$
X_0 = B_0\cup H_1\cup\cdots \cup H_n \cup \text{ (free handles)}.
$$

\smallskip\noindent
{\it \underline{Step 3}.}
The boundary of $X_0$ is an integral homology $3$-sphere $\Sigma$ with a free $\Z_2$-action.
Under the assumptions of \thmref{thm:EE}, Edmonds and Ewing proved that there exists a contractible $4$-manifold $Z$ with a locally linear $\Z_2$-action such that its boundary is $\Sigma$ with the given free $\Z_2$-action, and it has exact one fixed point.
Now $X=X_0\cup Z$ is the required manifold with the required action.
\medskip

Note that the action constructed above is smooth except near the final fixed point, that is, smooth on $X_0$.
We shall determine the sign assignment $\varepsilon$ on $X_0$. 
Note that each of $B_0$ and $H_1,\ldots, H_n$ has exact one fixed point.
We compare the value of $\varepsilon$ of the fixed point in each $H_i$ with that of $B_0$. 
It suffices to consider on $W=B_0\cup H_i$ for each $i$.

Recall that $W$ is given as $W=B_0\cup_{f_r}D^2\times D^2$, for a knot $K$ in $S_0=\partial B_0$ and a framing $r$.
Let $P$ be the fixed point in $B_0$ and $Q$ the fixed point in $D^2\times D^2$.
%
%
%
%
\begin{Proposition}\label{prop:framing}
Suppose $K$ is a trivial knot in $\partial B_0$ which bounds a $\Z_2$-invariant embedded disk $D_0$ in $B_0$ containing $P$.
If $r\equiv 2$ mod $4$, then $\varepsilon(P)=\varepsilon(Q)$. 
If $r\equiv 0$ mod $4$, then $\varepsilon(P)=-\varepsilon(Q)$. 
\end{Proposition}
\proof
Let $\iota\colon W\to W$ be the generator of the $\Z_2$-action.
In $H_i=D^2\times D^2$,  the disk $D_1:=D^2\times\{0\}$ is $\iota$-invariant and bounded by the knot $K$.
Let $N(D_0)$ be a tubular neighborhood of $D_0$. 
Then $N(D_0)$ can be identified with $D^2\times D^2$ via a $\Z_2$-equivariant map $\tilde{f}_0\colon D^2\times D^2\to N(D_0)$. 
The restriction of $\tilde{f}_0$ to $N(D_0)\cap \partial B_0$ gives a $0$-framing $f_0$ of $K$. 

Let $S$ be the $2$-sphere obtained by attaching $D_1$ to $D_0$ along $K$. 
Then a tubular neighborhood $N(S)$ of $S$ can be identified  with the manifold $D^2_0\times D^2\cup_{\bar{f}_r}D^2_1\times D^2$, where the attaching map $\bar{f}_r\colon \partial D_1\times D^2\to \partial D_0\times D^2$ is given by $\bar{f}_r (z,w) = (z,z^rw)$.   

The restriction of the frame bundle $F$ of $W$ to $S$  can be identified with the $\SO(4)$-bundle obtained by attaching $D_0\times \SO(4)$ with $D_1\times \SO(4)$ via the map $f^\prime\colon K=D_0\cap D_1\to \U(2)\subset \SO(4)$  given by 
$$
f^\prime (z) = 
\begin{pmatrix}
z^{-2} & 0 \\
0 & z^{-r}
\end{pmatrix},
$$
i.e., $F|_S = D_0\times \SO(4)\cup_{f^\prime_r}D_1\times \SO(4)$. 
In the above matrix, $z^{-2}$ comes from the Euler class of the tangent bundle of $S$, and $z^{-r}$ comes from the $r$-framing.
Note that the $\Z_2$-action on each $D_i\times \SO(4)$ is given by $\iota_*(z,v)=(-z,-v)$.

Take a path $\bar{s}$ in $S$ starting from $P$ and ending at $Q$. 
Then $\bar{s}$ and $\iota\bar{s}$ form a circle $\bar{C}$ in $S$, which divides $S$ into two disks: $S= D^+\cup D^-$. 
For $i=0,1$, let $s_i\colon \bar{s}\cap D_i\to F|_{D_i}=D_i\times \SO(4)$ be the constant section given by $s_i(t)=(t,I)$, where $I$ is the identity element of $\SO(4)$. 
Gluing $s_0$ and $s_1$, we obtain the lift $s$ of $\bar{s}$ to $F$.
Let $C=s\cup-\iota^* s$.

A spin structure on  $F|_{D^+}$ can be given by a trivialization of $F|_{D^+}$.
We shall concretely give such a trivialization as follows.
Let $D_i^+=D_i\cap D^+$ and $K^+=K\cap D^+$.
On $D_0^+$, assume the identity map as the trivialization $\phi_0 = \id \colon D_0^+\times \SO(4)\to F|_{D_0^+}=D^+_0\times \SO(4)$. 
Then, the trivialization $\phi_0$ on $D_0^+$ determines a trivialization on $K^+$ as $\phi^\prime\colon K^+\times \SO(4)\to F|_{K^+}=K^+\times \SO(4)$.
If we view $\phi^\prime$ {\it in} $D_1^+$, this is given by $\phi^\prime(z,v)=(z,f^\prime (z) v)$.
Then extend this trivialization $\phi^\prime$ over $D_1^+$ as $\phi_1\colon D_1^+ \times\SO(4)\to F|_{D_1^+}=D^+_1\times \SO(4)$. 
Gluing $\phi_0$ and $\phi_1$ along $K^+$, we obtain a trivialization $\phi\colon D^+\times \SO(4)\to F|_{D^+} = D_0^+\times \SO(4)\cup_{f^\prime_r}D_1^+\times \SO(4)$. 

Note that the trivialization $\phi$ above determines a double covering $\varphi\colon \hat{F}\to F$ on $D^+$ which gives the spin structure. 
Now we will apply \propref{prop:criterion}.
Let $p\colon D^+\times\SO(4)\to \SO(4)$ be the projection, 
and $C^\prime$ be the circle in $\SO(4)$ given by $C^\prime = p\circ\phi^{-1}(C)$.
Then $\varphi^{-1}(C)$ is connected if and only if $C^\prime$ represents the generator of $\pi_1(\SO(4))\cong\Z_2$.
Furthermore, by the construction above, we can see that $C^\prime$ represents the generator of $\pi_1(\SO(4))$ if and only if $(r+2)/2$ is odd.
\endproof

Thus, if each component of the link $L_T$ bounds an $\Z_2$-invariant embedded disk, the sign assignment $\varepsilon$ of $X_0$ is determined by diagonal entries $a_{ii}$ mod $4$ in the attaching matrix $(a_{ij})$.

%
%
\subsection{Nonsmoothable locally linear $\Z_2$-action on $K3\#K3$}\label{subsec:nonsmoothable}
%
%
Now, we construct a nonsmoothable locally linear $\Z_2$-action on $X=K3\#K3$.
Note that the intersection form of $X$ is isomorphic to $4E_8\oplus 6H$, where $H$ is the hyperbolic form.
By \thmref{thm:EE}, if we give an appropriate $\Z_2$-action on the intersection form, then we have a locally linear $\Z_2$-action on $X$.

Let $\Z_2$ act on $\Psi = 4E_8\oplus 6H$ as follows.
On a $2H$ summand, let $\Z_2$ act by permutation of two $H$'s. 
Similarly, on $2E_8\oplus 2E_8$ summand, let $\Z_2$ act by permutation of two $2E_8$'s, and on the rest $4H$ trivially.
The trivial part is denoted by $T$.

Let us consider the matrix,
\begin{equation}\label{eq:A}
A=
\begin{pmatrix}
0 & 1 & 1 & 1 & 1 & 1 & 1 & 1 \\
1 & 0 & 1 & 1 & 1 & 1 & 1 & 1 \\
1 & 1 & 0 & 1 & 1 & 1 & 1 & 1 \\
1 & 1 & 1 & 0 & 1 & 1 & 1 & 1 \\
1 & 1 & 1 & 1 & 2 & 1 & 1 & 1 \\
1 & 1 & 1 & 1 & 1 & 2 & 1 & 1 \\
1 & 1 & 1 & 1 & 1 & 1 & 2 & 1 \\
1 & 1 & 1 & 1 & 1 & 1 & 1 & 2 \\
\end{pmatrix}.
\end{equation}
Since the symmetric form represented by $A$ is even, indefinite and unimodular, it is isomorphic to $4H$.  
Therefore, we may assume $\Psi|_T$ is represented by the matrix $A$.
Furthermore, the matrix $A$ can be realized by a link whose every component bounds a $\Z_2$-invariant embedded disk as follows.
Let $p\colon S^3\to S^2$ be the Hopf fibration. Take distinct $8$ points in $S^2$. Then the inverse image of $8$ points by $p$ forms a required link.

As in \subsecref{subsec:EE}, by equivariant handle construction, we can construct a smooth action on a manifold $X_0$ with a boundary, and \thmref{thm:EE} says that this action is extended to the whole $X$ as a locally linear action.

If the smooth action on $X_0$ is extended {\it smoothly} on the whole $X$, then we have $\sum\varepsilon(P)=\pm 8$ by \propref{prop:K3K3}.
However, this is a contradiction because $\sum\varepsilon(P)$ can not be $\pm 8$ by \propref{prop:framing} for the action constructed on the matrix $A$.
Thus, the smooth action on $X_0$ can not be extended to the whole $X$ {\it smoothly}.

Now, we claim more strongly that the locally linear action constructed above is {\it nonsmoothable}. 
First, we shall clarify what nonsmoothable means.
Let $X$ be an oriented {\it topological} manifold and $G$ a finite group. 
If $X$ admits  a {\it smooth} structure and a smooth structure $\sigma$ is specified, then we write the manifold with the smooth structure $\sigma$ by $X_\sigma$.
Let $LL(G,X)$ be the set of equivalence classes of orientation-preserving locally linear $G$-actions on $X$. 
Here, two locally linear actions are said equivalent if there exists a {\it homeomorphism} $f$ of $X$ such that one action is conjugate to the other by $f$.
Similarly, let $C^{\infty}(G,X_{\sigma})$ be the set of equivalence classes of orientation-preserving smooth $G$-actions on $X_{\sigma}$. 
Here, the equivalence of two smooth actions is given via a {\it diffeomorphism} of $X_{\sigma}$.
Then we have a forgetful map $\Phi_{\sigma}\colon C^{\infty}(G,X_{\sigma})\to LL(G,X)$ forgetting the smooth structure.
\begin{Definition}
A locally linear action is called {\it nonsmoothable} with respect to the smooth structure $\sigma$ if its class is not contained in the image of $\Phi_{\sigma}$.
\end{Definition}
Let $X_{\sigma}$ be a closed  smooth $4$-manifold and $G=\Z_2$.
Since a spin structure on $X_{\sigma}$ is defined from the tangent bundle, the definition implicitly uses the smooth structure on $X_{\sigma}$. 
Hence, at present, the sign assignment $\varepsilon$ seems to  depend on classes of $C^{\infty}(G,X_{\sigma})$. 

We claim that the sign assignment  can be defined for locally linear actions, and it depends only on classes in $LL(G,X)$.
As in \cite{E1}, topological spin structures on topological manifolds can be defined as follows.
An oriented topological $n$-manifold $M$ has the tangent micro bundle $\tau M$. 
By Kister-Mazur's theorem~\cite{Kister},  $\tau M$ can be identified with a $\R^n$-bundle $\tau^\prime M$ whose structural group is the topological group $\STop(n)$ which consists of orientation-preserving homeomorphisms of $\R^n$ preserving the origin. 
It is known that, when $n\geq 2$, the inclusion $\SO(n)\to\STop(n)$ induces isomorphisms of $\pi_0$, $\pi_1$ and $\pi_2$. (See \cite{KS} and \cite{FQ}.)
Hence there is the unique double covering $\varphi_0\colon\SpinTop(n)\to\STop(n)$.
Let $F$ be the principal $\STop(n)$-bundle of ``frames'' of $\tau^\prime M$.
Then, a topological spin structure is given by a double covering $\varphi \colon \hat{F}\to F$ whose restriction to each fibre is $\varphi_0$.

Suppose that we are given a locally linear $\Z_2$-action on a simply-connected topological $4$-manifold $X$ with isolated fixed points.
Then the tangent micro bundle $\tau X$ is $\Z_2$-locally linear in the sense of \cite{LR}, and the corresponding $\R^n$-bundle $\tau^\prime X$ becomes a $\Z_2$-equivariant bundle.
Let $F$ be  the principal $\STop(4)$-bundle of frames of $\tau^\prime X$.
Then, a $\Z_2$-action on $F$ is induced from that on $\tau^\prime X$.
By considering as in \cite{AB} and \cite{E1}, we see that the $\Z_2$-action is even type if and only if fixed points are isolated.

Now, we can define the sign assignment by using Atiyah-Bott's criterion itself  on the topological spin structure $\varphi \colon\hat{F}\to F$ as follows.
Since the action is assumed locally linear, on the fiber of $F$ over each fixed point $P$, there is a point $y_P$ which is mapped to $-y_P$ by the $\Z_2$-action, where ``$-$'' means the multiplication of $-1$ which is considered as an element of $\SO(4)\subset \STop(4)$. 
For distinct fixed points $P$ and $Q$, by taking a path $s$ connecting such a $y_P$ with such a $y_Q$, we can define the sign assignment as in \propref{prop:criterion}:
\begin{Definition} 
For each pair $(P,Q)$ of fixed points, let $s$ be a path in $F$ as above, and $C$ the circle formed by $s$ and $-\iota_*s$. 
Define $\varepsilon^\prime (P,Q)$ by 
\begin{align*}
\varepsilon^\prime (P,Q) &= 1,\quad \text{if $\varphi^{-1}(C)$ has $2$ components,}\\
\varepsilon^\prime (P,Q) &= -1,\quad \text{if $\varphi^{-1}(C)$ is connected.}
\end{align*}
\end{Definition}
This definition does not depend on smooth structures, and is well-defined if $X$ is simply-connected.
Furthermore, if the action is realized by a smooth action, 
\begin{equation}\label{eq:ep}
\varepsilon^\prime (P,Q)=\varepsilon(P)\varepsilon(Q).
\end{equation}

Now, we complete the proof of \thmref{thm:nonsmoothable}.
\proof[Proof of \thmref{thm:nonsmoothable}]
We have a locally linear action constructed from the matrix $A$ in \eqref{eq:A}.
If this action is smoothable for a smooth structure, then $\sum\varepsilon(p)$ must be $\pm 8$ by \propref{prop:K3K3}.
However, this is impossible by \propref{prop:framing} and the relation \eqref{eq:ep}.
\endproof

%
%
\subsection{Nonsmoothable locally linear $\Z_2$-action on $K3$}\label{subsec:K3}
%
%
The argument in previous subsections enables us to prove \thmref{thm:K3}. 
The proof is parallel to that of \thmref{thm:nonsmoothable}.

Since both $k_+$ and $k_-$ are even and $\Sign(K3)=-16$, \eqref{eq:G-ind} implies the following constraint.
\begin{Proposition}\label{prop:K3}
Let $X$ be a simply-connected spin manifold which has the same intersection form as $K3$.
Suppose that an orientation-preserving smooth $\Z_2$-action on $X$ is given, and let $c_0$ be the $\Spinc$-structure which is determined by a spin structure.
If $b_+^{\Z_2}=3$ and fixed points are discrete, then the lift of $\Z_2$-action to $c_0$ satisfies, 
\begin{equation*}
\sum_{p} \varepsilon(p)=\pm  8.
\end{equation*}
\end{Proposition}

Now, we shall construct a nonsmoothable locally linear action on $K3$.
The intersection form of $X=K3$ is isomorphic to $2E_8\oplus 3H$. 
Let $\Z_2$ act on the $2E_8$ summand by permutation, and on the $3H$ summand trivially.
Let us consider the matrix,
\begin{equation*}
B=
\begin{pmatrix}
0 & 1 & 1 & 1 & 1 & 1\\
1 & 0 & 1 & 1 & 1 & 1\\
1 & 1 & 0 & 1 & 1 & 1\\
1 & 1 & 1 & 2 & 1 & 1\\
1 & 1 & 1 & 1 & 2 & 1\\
1 & 1 & 1 & 1 & 1 & 2 
\end{pmatrix}.
\end{equation*}
Since the symmetric form represented by $B$ is even, indefinite and unimodular, $B$ is isomorphic to $3H$. 
By using the Hopf fibration, we can find  a link  representing the matrix $B$ whose every component bounds a $\Z_2$-invariant embedded disk.
As in \subsecref{subsec:EE}, by equivariant handle construction, we obtain a locally linear $\Z_2$-action associated to $B$.
If this action is smoothable, the lift of the action to the spin structure should satisfy $\sum\varepsilon(p)=\pm 8$ by \propref{prop:K3}.
However this is impossible by \propref{prop:framing} and the relation \eqref{eq:ep}.
Thus, \thmref{thm:K3} is established.


\begin{thebibliography}{99}

\bibitem{AB} M.~F.~Atiyah and R.~Bott, 
{\it A Lefschetz fixed point formula for elliptic complexes \Romnum{2} Applications}, 
Ann. Math. {\bf 88} (1968), 451--491. 

\bibitem{B0} S.~Bauer,
{\it A stable cohomotopy refinement of Seiberg-Witten invariants. \Romnum{2}}, 
Invent. Math. {\bf 155} (2004), no. 1, 21--40. 

\bibitem{B1} S.~Bauer,
{\it Refined Seiberg-Witten invariants},
in: S.~K.~Donaldson, Y.~Eliashberg, M. Gromov (edits.) {\it Different faces of Geometry}, Internet. Math. Ser., Kluwer Academic/Plenum Publishers, New York, 2004.

\bibitem{B2} S.~Bauer,
{\it Almost complex $4$-manifolds with vanishing first Chern class},
preprint, math.GT/0607714.

\bibitem{BF} S.~Bauer and M.~Furuta, 
{\it A stable cohomotopy refinement of Seiberg-Witten invariants: \Romnum{1}},
Invent. Math. {\bf 155} (2004), 1--19.

\bibitem{Bredon} G.~E.~Bredon,
{Equivariant cohomology theories},
Lecture Notes in Math., vol. 34, Springer, Berlin, 1967. 

\bibitem{Bryan} J.~Bryan,
{\it Seiberg-Witten theory and $\Z/2^p$ actions on spin $4$-manifolds},
Math. Res. Lett. {\bf 5} (1998), no. 1-2, 165--183. 

\bibitem{Dieck} T.~tom Dieck,
{Transformation groups},
de Gruyter Stud. Math., vol. 8, Springer, Berlin, 1987.

\bibitem{E1} A.~L.~Edmonds,
{\it Aspects of group actions on four-manifolds},
Topology Appl. {\bf 31} (1989), no. 2, 109--124. 

\bibitem{EE} A.~L.~Edmonds and J.~H.~Ewing,
{\it Realizing forms and fixed point data in dimension four},
Amer. J. Math. {\bf 114} (1992), 1103--1126.

\bibitem{Fang} F.~Fang,
{\it Smooth group actions on $4$-manifolds and Seiberg-Witten invariants},
International J.~Math. {\bf 9}, No.8 (1998) 957--973. 

\bibitem{FQ} M.~H.~Freedman and F.~Quinn,
Topology of 4-manifolds,
Princeton Mathematical Series, 39. Princeton University Press, Princeton, NJ, 1990.

\bibitem{FKM} M.~Furuta, Y.~Kametani, N.~Minami,
{\it Stable-homotopy Seiberg-Witten Invariants for Rational Cohomology $K3\#K3$'s},
Journal of Mathematical Sciences, The University of Tokyo {\bf 8} (2001), No. 1, 157-176.

\bibitem{HY} A.~Hattori and T.~Yoshida,
{\it Lifting compact group actions in fiber bundles}, 
Japan. J. Math. (N.S.) {\bf 2} (1976), no. 1, 13--25. 

\bibitem{Hu} S.~Hu,
Homotopy theory, 
Pure and Applied Mathematics, Vol. VIII, Academic Press, New York, 1959.

\bibitem{KS} R.~C.~Kirby and L.~C.~Siebenmann, 
Foundational essays on topological manifolds, smoothings, and triangulations. 
Annals of Mathematics Studies, No. 88. Princeton University Press, Princeton, N.J.; University of Tokyo Press, Tokyo, 1977. 

\bibitem{Kister} J.~M.~Kister,
{\it Microbundles are fibre bundles}.
Ann. of Math. (2) {\bf 80} (1964) 190--199.

\bibitem{LR} R.~Lashof and M.~Rothenberg,
{\it $G$-smoothing theory}, 
Algebraic and geometric topology, Proc. Sympos. Pure Math. XXXII, Part 1, pp. 211--266, Amer. Math. Soc., Providence, R.I., 1978. 

\bibitem{LN} X.~Liu and N.~Nakamura, 
{\it Pseudofree $\Z/3$-actions on $K3$ surfaces}, 
Proc. Amer. Math. Soc. {\bf 135} (2007), no. 3, 903--910. 

\bibitem{Mac}S.~Mac Lane, 
Homology, Springer-Verlag, New York, 1975.

\bibitem{LN2} X.~Liu and N.~Nakamura, 
{\it Nonsmoothable group actions on elliptic surfaces}, 
Topology Appl. {\bf 155} (2008), 946--964.

\bibitem{M2} T.~Matumoto,
{\it Equivariant cohomology theories on $G$-CW complexes},
Osaka J. Math. {\bf 10} (1973), 51--68.

\bibitem{May} J.~P.~May et al., 
Equivariant homotopy and cohomology theory, 
CBMS, Vol. 91, American Mathematical Society, Providence, RI, 1996

\bibitem{Morrison} D.~R.~Morrison,
{\it On $K3$ surfaces with large Picard number},
Invent. Math. {\bf 75} (1984), no. 1, 105--121. 

\bibitem{Np} N.~Nakamura,
{\it Mod $p$ vanishing theorem of Seiberg-Witten invariants for $4$-manifolds with $\Z_p$-actions},
Asian J. Math. {\bf 9}, (2006), no.4, 731--748.

\bibitem{Sz} M.~Szymik,
{\it Bauer-Furuta invariants and Galois symmetries},
preprint.

\end{thebibliography}
\end{document}